\documentclass[a4paper]{article}
\usepackage{amssymb, latexsym, amsmath, amsthm}
\usepackage{a4wide}

\newtheorem{theorem}{Theorem}[section]
\newtheorem{lemma}[theorem]{Lemma}
\newtheorem{proposition}[theorem]{Proposition}

\newtheorem{corollary}[theorem]{Corollary}

\theoremstyle{definition}

\newcommand{\K}{\ensuremath{\mathbb{K}}}

\newcommand{\Z}{\ensuremath{\mathbb{Z}}}
\newcommand{\N}{\ensuremath{\mathbb{N}}}

\newcommand{\q}{\quad}
\newcommand{\qq}{\qquad}

\newcommand{\germ}{\mathfrak}

\newcommand{\sln}{\mathfrak{sl}_{n}}
\newcommand{\slf}{\mathfrak{sl}^{+}_{4}}

\newcommand{\U}{U_{q}({\slf})}

\newcommand{\aut}{{\rm Aut}_{\K}\U}
\newcommand{\qhat}{\hat{q}}

\newcommand{\fract}{\mathrm{Frac}}
\newcommand{\der}{\mathrm{Der}}

\def\qtor{P(\Lambda)}

\begin{document}
\title{On the Hochschild cohomology and the automorphism group of $\U$}
\author{St\'ephane Launois\thanks{This research was supported by a  Marie Curie Intra-European
  Fellowship within the $6^{\mbox{\tiny th}}$ European Community  Framework
  Programme}\q and\q Samuel A. Lopes\thanks{Work partially supported  by \emph{Centro de Matem\'atica da Universidade do Porto} (CMUP),  financed by FCT (Portugal) through the programmes POCTI (\emph {Programa Operacional Ci\^encia, Tecnologia, Inova\c c\~ao}) and POSI  (\emph{Programa Operacional Sociedade da Informa\c c\~ao}), with  national and European Community structural funds.}}
\date{}
\maketitle

\begin{abstract}
We compute the automorphism group of the $q$-enveloping algebra $\U$  of the
nilpotent Lie algebra of strictly upper triangular matrices of size
$4$. The result obtained gives a positive answer to a conjecture of
Andruskiewitsch and Dumas. We also compute the derivations of this
algebra and then show that the Hochschild cohomology group of degree 1
of this algebra is a free (left)-module of rank 3 (which is the rank
of the Lie algebra $\germ{sl}_4$) over the center of $\U$.
\end{abstract}

\noindent {\it Keywords:} quantized enveloping algebra; automorphisms;
derivations; Hochschild cohomology.
\\$ $
\\{\it 2000 Mathematics Subject Classification:} 17B37, 16W20, 16W25,  16E40.

\section*{Introduction}

Let $\K$ be a field, $\mathcal{L}$ a Lie algebra over the  $\K$ and $U(\mathcal{L})$ its enveloping algebra. The group ${\rm Aut}_{\K}U(\mathcal{L})$ of $\K$-algebra automorphisms of $U(\mathcal{L})$ is still for the most part unknown (except in particular instances, e.g.\ $\dim\mathcal{L}\leq 2$). For example, if 
$\mathcal{L}$ is the two-dimensional abelian Lie algebra, then $U(\mathcal{L})$ is the polynomial algebra in two indeterminates $x_{1}$, $x_{2}$, whose group of automorphisms is generated by the \emph{elementary} automorphisms of the form
\begin{equation*}
x_{i}\mapsto \lambda x_{i}+f(x_{j}), \q x_{j}\mapsto x_{j} \q\q\q (i\neq j)
\end{equation*}
with $\lambda\in\K^*$ and $f(x_{j})$ a polynomial in the variable $x_{j}$ (\cite{hweJ42}, \cite{wVDK53}). In contrast with this simple description, the conjecture that the polynomial algebra in three variables over 
$\K$ has \emph{wild} automorphisms (i.e.\ automorphisms not of the above type) has recently been settled (see~\cite{SU04}) assuming $\K$ has characteristic $0$. Another example is the enveloping algebra of $\mathfrak{sl}_{2}$, which is known to have wild automorphisms by a result of Joseph~\cite{aJ76}.

Pertaining more to what is studied in this paper is the enveloping algebra of the three-dimensional Heisenberg Lie algebra, which is given by generators $x$, $y$ and $z$, subject to the relations
\begin{equation*}
[x, y]=z,\q [z, x]=0=[z, y].
\end{equation*}
This algebra can also be seen as the enveloping algebra of the Lie algebra $\mathfrak{sl}^{+}_{3}$ of strictly upper triangular matrices of size $3$. The infinite dimensional simple quotients of $U(\mathfrak{sl}^{+}_{3})$ are isomorphic to the first Weyl algebra $\mathbb{A}_{1}(\K)$, whose group of automorphisms was described by Dixmier in~\cite{jD68}. Yet, the full group of automorphisms of $U(\mathfrak{sl}^{+}_{3})$ remains to be described, and Alev~\cite{jA86} proved the existence of wild automorphisms of this algebra.

Unlike the classical scenario, quantum algebras are believed to possess less symmetry (see~\cite[1.1]{kG03}) and the group of automorphisms of several algebras of this kind has been computed successfully. Making use of a general result relating automorphisms  and  derivations of $\N$-graded algebras, Alev and Chamarie~\cite{AC92} described the automorphism group of quantum affine space, of the algebra of $2\times 2$ quantum matrices and of the quantized enveloping algebra $U_{q}(\mathfrak{sl}_{2})$. Also, in~\cite{AD96} the authors found the automorphism groups of the quantum Weyl algebra, the Weyl-Hayashi algebra, the quantum Heisenberg algebra $U_{q}(\mathfrak{sl}^{+}_{3})$ (see also~\cite{pC95}) and of other related algebras. Here the methods used included describing the set of normal elements of the algebras involved and using appropriate filtrations to carry out computations. In~\cite{lR96}, Rigal used the invariance under automorphisms of the set of height $1$ prime ide
 als
  of quantum Weyl algebras to describe their automorphism group. Related methods were employed by  G\'omez-Torrecillas and Kaoutit~\cite{GTEK02} regarding the coordinate ring of quantum symplectic space, and by Lenagan and the first author~\cite{LnLe05pr} regarding the algebra of non-square quantum matrices. In all of these cases, the automorphism group of the algebras involved does not differ from the natural torus which acts diagonally on the generators by more than a finite group and perhaps a copy of $\Z$.

In their paper~\cite{ADpr}, Andruskiewitsch and Dumas conjectured that, given a finite-dimensional complex simple Lie algebra $\mathfrak{g}$ with triangular decomposition 
$\mathfrak{g}=\mathfrak{g}^{-}\oplus\mathfrak{h}\oplus\mathfrak{g}^{+}$, then 
${\rm Aut}_{\K}U_{q}(\mathfrak{g}^{+})$, the group of $\K$-algebra automorphisms of the quantized enveloping algebra of the nilpotent Lie algebra $\mathfrak{g}^{+}$, is isomorphic to the semi-direct product of the torus $(\K^{*})^n$ ($n$ being the rank of $\mathfrak{g}$) with the group of order $1$, $2$ or $3$ generated by the diagram automorphism of $\mathfrak{g}^{+}$, see~\cite[Prob.\ 1]{ADpr}. This conjecture holds for 
$\mathfrak{g}^{+}=\mathfrak{sl}^{+}_{3}$ (\cite{pC95}, \cite{AD96}) and  recently the first author proved in~\cite{sLnpr} that it holds as well in the $B_{2}$ case, i.e., with $\mathfrak{g}^{+}=\mathfrak{so}^{+}_{5}$.

In this paper we settle the conjecture of Andruskiewitsch and Dumas in the $A_{3}$ case, so that $\mathfrak{g}^{+}=\mathfrak{sl}^{+}_{4}$ is the Lie algebra of strictly upper triangular matrices of size $4$. We also compute the Lie algebra of derivations and the first Hochschild cohomology group of  $U_{q}(\mathfrak{sl}^{+}_{4})$, which is shown to be a free module of rank $3$ over the center of $U_{q}(\mathfrak{sl}^{+}_{4})$.

Let us briefly summarise  what is done in the paper. There exist normal elements $\Delta_{1}$, $\Delta_{2}$ and $\Delta_{3}$ such that the center of $U_{q}(\mathfrak{sl}^{+}_{4})$ is the polynomial algebra in the variables $z_{1}=\Delta_{1}\Delta_{3}$ and $z_{2}=\Delta_{2}$. Given an automorphism $\phi$ of $U_{q}(\mathfrak{sl}^{+}_{4})$, our strategy is to show that, up to the diagram automorphism and the diagonal action of the torus $(\K^{*})^{3}$ on the Chevalley generators of $U_{q}(\mathfrak{sl}^{+}_{4})$, $\phi$ fixes $\Delta_{1}$, $\Delta_{2}$ and $\Delta_{3}$. Then, by using degree arguments, we conclude that $\phi$ is the identity.

The difficulty that arises is in showing that the central element $\Delta_{2}$ is fixed. Hence we use the methods of~\cite{AC92} and~\cite{LL06} and determine the derivations of $U_{q}(\mathfrak{sl}^{+}_{4})$. To do this, we first apply the deleting derivations algorithm of Cauchon~\cite{gC03} so that, after suitably localising, we can embed $U_{q}(\mathfrak{sl}^{+}_{4})$ in a quantum torus $\qtor$. Extending a derivation $D$ of $U_{q}(\mathfrak{sl}^{+}_{4})$ to $\qtor$ we obtain, by a result of Osborn and Passman~\cite{OP95}, a decomposition
\begin{equation*}
D=\mathrm{ad}_{x}+\theta
\end{equation*}
with $x\in\qtor$ and $\theta$ a central derivation of $\qtor$. Using a sort of restoring derivations algorithm, we finish by deducing that $x\in U_{q}(\mathfrak{sl}^{+}_{4})$ and that $\theta$ sends each Chevalley generator of $U_{q}(\mathfrak{sl}^{+}_{4})$ to a multiple of itself by a central element of $U_{q}(\mathfrak{sl}^{+}_{4})$.

$ $

\begin{flushleft}
\textbf{Acknowledgments.} Most of this work was done during a visit of the second author to the School of Mathematics at the University of Edinburgh. He would like to express his gratitude for the hospitality received, especially from T.H.\ Lenagan, L.\ Richard and the first author.
\end{flushleft}

\section{Basic aspects of $\U$}\label{S:int}

Let $\K$ be a field of characteristic $0$ and fix a parameter $q\in\K^ {*}$ which we assume is not a root of unity. Consider, for $n\geq 2$,  the Lie algebra $\germ{sl}_{n}$ of $n\times n$ matrices of trace $0$  and its maximal nilpotent subalgebra $\germ{sl}_{n}^+$ consisting of  the strictly upper triangular matrices of size $n$.

Throughout this paper $\N$ is the set of nonnegative integers. For $k \in\N$, the $q$-integer $[k]$ is defined by $[k]=\frac{q^{k}-q^{-k}} {q-q^{-1}}$ and we use the notation $\widehat{q}=q-q^{-1}$.


\subsection{$q$-Serre relations}\label{SS:int:qS}

The algebra $\U$ is the $q$-deformation of the universal enveloping  algebra of the nilpotent Lie algebra $\slf$. It is the unital  associative $\K$-algebra with generators $e_{1}$, $e_{2}$ and $e_{3} $, subject to the quantum Serre relations:
\begin{align}
e_{1}e_{3}-e_{3}e_{1}&=0\\
e_{i}^2 e_{j}-(q+q^{-1})e_{i}e_{j}e_{i}+e_{j}e_{i}^2&=0\qq \mbox{if\  \ $|i-j|=1$.}
\end{align}


\subsection{Weight space decomposition}\label{SS:int:wsd}

Let $Q=\Z^3$ be the free abelian group with canonical basis $\alpha_ {1}$, $\alpha_{2}$, $\alpha_{3}$ and $Q^{+}=\N^{3}$ be its submonoid.  Since the quantum Serre relations are homogeneous in the given  generators, there is a $Q^+$-grading on $\U$ obtained by assigning to  $e_{i}$ degree $\alpha_{i}$. We use the terminology \emph{weight}  instead of degree for this grading, and write $wt (u)=\beta$ if $u\in \U$ has weight $\beta$.


\subsection{PBW basis}\label{SS:int:pbw}

Several authors have constructed PBW bases for quantized enveloping  algebras (e.g.~\cite{hY89}, \cite{mT90}, \cite{cR96}). It will be  convenient for us to use the following construction:
$$
\begin{array}{l@{\qq\q}l}
X_{1}=e_{1},&X_{2}=e_{1}e_{2}-q^{-1}e_{2}e_{1},\\[5pt]
X_{4}=e_{2},&X_{5}=e_{2}e_{3}-q^{-1}e_{3}e_{2},\\[5pt]
X_{6}=e_{3},&X_{3}=e_{1}X_{5}-q^{-1}X_{5}e_{1}.
\end{array}
$$
Then, the set of monomials
$\displaystyle \left\{ X^{b_{1}}_{1}\cdots X^{b_{6}}_{6}\mid b_{i}\in \N\right\}$ is a linear basis of $\U$. Notice that all $X_{i}$ are  weight vectors.


\subsection{Ring theoretical properties of $\U$}\label{SS:int:rtp}

It was seen in~\cite{cR96} (see also~\cite[I.6.10]{BG02} and  references therein) that $\U$ is an iterated skew polynomial ring. In  terms of the PBW basis described above, we have
\begin{equation}\U=\K[X_{1}][X_{2};\tau_{2}][X_{3};\tau_{3}][X_{4}; \tau_{4}, \delta_{4}][X_{5};\tau_{5}, \delta_{5}][X_{6};\tau_{6},  \delta_{6}],\end{equation}
with $\tau_{i}$ a $\K$-algebra automorphism and $\delta_{i}$ a $\K$- linear (left) $\tau_{i}$-derivation of the appropriate subalgebra.  Thus $\U$ is a Noetherian domain.

So that we can easily compute in $\U$, and also because this  information will be needed in Section~\ref{S:der}, we specify these  automorphisms and skew-derivations below by giving their values on the $X_ {j}$ ($\delta_{i}(X_ {j})=0$ unless otherwise specified):
\begin{gather*}
 \tau_{2}(X_{1})=q^{-1}X_{1}
\end{gather*}
\begin{gather*}
   \tau_{3}(X_{1})=q^{-1}X_ {1},     \quad\quad   \tau_{3}(X_{2})=q^{-1}X_{2}    
\end{gather*}
\begin{gather*}
   \tau_{4}(X_{1})=qX_{1},     \quad\quad     \tau_{4}(X_{2})=q^{-1}X_{2},       \quad\quad   \tau_{4}(X_{3})=X_{3},  \quad\quad   \delta_{4}(X_{1})=-qX_{2}
\end{gather*}
\begin{gather*}
\tau_{5}(X_{1})=qX_{1},  \quad  \quad    \tau_{5}(X_{2})=X_{2},  \quad   \quad  \tau_{5}(X_{3})=q^{-1}X_{3},  \\ 
      \tau_{5}(X_{4})=q^{-1}X_{4},  \quad \quad   \delta_{5}(X_{1})=-qX_{3},  \quad \quad \delta_{5}(X_ {2})=-\qhat X_{3}X_{4}
\end{gather*}
\begin{gather*}
\tau_{6}(X_{1})=X_{1},  \quad \quad   \tau_{6}(X_{2})=qX_{2},  \quad \quad  \tau_{6}(X_{3})=q^{-1}X_{3}\\
 \tau_{6}(X_{4})=qX_{4},  \quad \quad   \tau_{6}(X_{5})=q^{-1}X_{5},  \quad \quad \delta_{6}(X_{2})=-qX_{3},  \quad \quad \delta_{6}(X_{4})=-qX_{5}
\end{gather*}
Furthermore, for $4\leq i\leq 6$, $\tau_{i}\circ\delta_{i}=q^{-2} \delta_{i}\circ\tau_{i}$, so the theory of deleting derivations of~ \cite{gC03} applies to $\U$. In particular, as shown in~\cite{cR96},  all prime ideals of $\U$ are completely prime.


\subsection{Normal elements and the center}\label{SS:int:nor}

The elements $a, b\in\U$ are said to \emph{$q$-commute} if there is  an integer $\lambda$ such that $ab=q^{\lambda}ba$. If $u$ $q$-commutes with the generators $e_{i}$ of $\U$ then we say that $u$ is  \emph{$q$-central}. Clearly, $q$-central elements are normal and  Caldero~\cite[Prop.~2.1]{pC95} has shown the reciprocal of this  statement, so that the normal elements of $\U$ are just the $q$-central ones.

The following theorem was established (in the more general context of  $U_{q}(\sln^{+})$) independently by Alev and Dumas~\cite{AD94} and by  Caldero~\cite{pC94a, pC95}.

\begin{theorem}\label{T:int:nor}
There exist $q$-central weight elements $\Delta_{i}\in\U$, $i=1, 2, 3 $, such that:
\begin{itemize}\item[\textup{(a)}] $\Delta_{2}$ is central and
\begin{itemize}\item[\textup{(i)}] $e_{2}$ commutes with $\Delta_{i} $, for all $i=1, 2, 3$;
\item[\textup{(ii)}] $e_{1}\Delta_{1}=q\Delta_{1}e_{1}$, $e_{1}\Delta_ {3}=q^{-1}\Delta_{3}e_{1}$;
\item[\textup{(iii)}] $e_{3}\Delta_{1}=q^{-1}\Delta_{1}e_{3}$, $e_{3} \Delta_{3}=q\Delta_{3}e_{3}$;\end{itemize}
\item[\textup{(b)}] The subalgebra $\K[\Delta_{1}, \Delta_{2}, \Delta_ {3}]$ generated by the $\Delta_{i}$ is a (commutative) polynomial  algebra  in $3$ variables.
\item[\textup{(c)}] The center $Z_{q}(\slf)$ of $\U$ is the  polynomial algebra in the variables $z_{1}=\Delta_{1}\Delta_{3}$ and  $z_{2}=\Delta_{2}$.\end{itemize}
\end{theorem}

The set of $q$-central elements of $\U$ was also described by Caldero  (see for example~\cite[Th\'e.~2.2]{pC95}) in terms of the $\Delta_{i} $ and the longest element of the Weyl group of $\mathfrak{sl}_{4}$  (in the notation of~\cite{pC95}, $\Delta_{i}=e_{s(\varpi_{4-i})}$).  It follows from his analysis that every $q$-central element is an  element of $\K[\Delta_{1}, \Delta_{2}, \Delta_{3}]$. So let $p=\sum_ {j}c_{j}\Theta_{j}$ be $q$-central, with each $c_{j}\in\K^{*}$ and   the $\Theta_{j}$ distinct monomials in the $\Delta_{i}$. Take $\lambda \in\Z$ so that $e_{1}p=q^{\lambda}pe_{1}$. By Theorem~\ref{T:int:nor} (a), each $\Theta_{j}$ is $q$-central, so it must be that $e_{1} \Theta_{j}=q^{\lambda}\Theta_{j}e_{1}$ for all $j$, as $\U$ is a  domain and the $\Theta_{j}$ are distinct. Assume $\lambda\geq 0$ and  write $\Theta_{j}=\Delta_{1}^{\alpha}\Delta_{2}^{\beta}\Delta_{3}^ {\gamma}$. Then, once more by Theorem~\ref{T:int:nor}(a), $\lambda= \alpha-\gamma$ and so $\Theta!
 _{j}=\Delta_{1}^{\lambda}u_{j}$ with $u_ {j}=z_{1}^{\gamma}z_{2}^{\beta}$ central. Since $j$ was arbitrary, we  deduce that $p$ is the product of $\Delta_{1}^{\lambda}$ and a  central element. Had we assumed $\lambda\leq 0$, we would have  obtained an analogous statement with $\Delta_{1}^{\lambda}$ replaced  by $\Delta_{3}^{-\lambda}$. Conversely, it is clear that all elements  of $\Delta_{i}^{c}Z_{q}(\slf)$ are $q$-central, for $c\in\N$ and $i\in \{ 1, 3 \}$, so we have established the following:

\begin{lemma}\label{L:int:nor}
Let $u\in\U$ be normal. Then there exists a central element $z$, a  nonnegative integer $c$ and $i\in\{1, 3\}$ such that $u=\Delta_{i}^{c} z$.
\end{lemma}

In terms of the PBW basis we are using, the $\Delta_{i}$ are given by  the formulae (see~\cite[Sec.\ 4]{pC94a} or~\cite[Sec.\ 4.1]{sLpr05}  but notice that we have ordered the PBW basis elements differently):
\begin{align}\label{E:int:nor:d1}
\Delta_{1}&=X_{3},\\\label{E:int:nor:d2}
\Delta_{2}&=X_{2}X_{5}-qX_{3}X_{4},\\\label{E:int:nor:d3}
\Delta_{3}&=\qhat{\,}^{2}X_{1}X_{4}X_{6}-q\qhat X_{2}X_{6} -q\qhat X_ {1}X_{5}+q^{2}X_{3}.
\end{align}


\section{The automorphism group of $\U$}\label{S:aut}

In this section we compute the group of algebra automorphisms of $\U$  and confirm the conjecture of Andruskiewitsch and Dumas~\cite{ADpr}  for this case. Let $\aut$ denote this group. We shall show that $\aut $ is the semi-direct product of the $3$-torus $(\K^{*})^{3}$ and the  group of order two generated by the diagram automorphism of $\U$.

Let $\mathcal{H}=(\K^{*})^{3}$. Each $\bar{\lambda}=(\lambda_{1},  \lambda_{2}, \lambda_{3})\in\mathcal{H}$ determines an algebra  automorphism $\phi_{\bar{\lambda}}$ of $\U$  with $\phi_{\bar {\lambda}}(e_{i})=\lambda_{i}e_{i}$ for $i=1, 2, 3$, with inverse
$\phi^{-1}_{\bar{\lambda}}=\phi_{\bar{\lambda}^{-1}}$. Hence we think  of $\mathcal{H}$ as a subgroup of $\aut$ via this correspondence.  There is also a diagram automorphism $\eta$ of $\U$ arising from the  symmetry of the Dynkin diagram of type $A$, and defined on the  generators by $\eta (e_{i})=e_{4-i}$. Notice that $\eta^{2}$ is the  identity morphism and that, up to nonzero scalars, $\eta$ permutes $ \Delta_{1}$ and $\Delta_{3}$, and fixes $\Delta_{2}$. Finally, as is  to be expected,
\begin{equation}\label{E:aut:phieta}\eta\circ\phi_{(\lambda_{1},  \lambda_{2}, \lambda_{3})}\circ\eta^{-1}=\phi_{(\lambda_{3}, \lambda_ {2}, \lambda_{1})}.\end{equation}

\subsection{An $\N$-grading on $\U$}\label{S:aut:ngr}

In addition to the weight space decomposition of Section~\ref {SS:int:wsd}, $\U$ has an $\N$-grading induced by the monoid  homomorphism $a\alpha_{1}+b\alpha_{2}+c\alpha_{3}\mapsto a+b+c$, from  $Q^{+}$ to $\N$. Let
\begin{equation}\U=\bigoplus_{i\in\N}U_{i}\end{equation}
be the corresponding decomposition, with $U_{i}$ the subspace of  homogeneous elements of degree $i$. In particular, $U_{0}=\K$ and $U_ {1}$ is the $3$-dimensional space spanned by the generators $e_{1}, e_ {2}, e_{3}$. For $t\in\N$ set $U_{\geq t}=\bigoplus_{i\geq t}U_{i}$  and define   $U_{\leq t}$ similarly.

We say that the nonzero element $u\in\U$ has degree $t$, and write $ \mathrm{deg} (u)=t$, if $u\in U_{\leq t}\setminus U_{\leq t-1}$  (using the convention that $U_{\leq -1}=\{ 0 \}$). In such a case, if  $u=\sum_{0\leq i\leq t} u_{i}$ with $u_{i}\in U_{i}$ and $u_{t}\neq 0 $, we set $\bar{u}=u_{t}$. By definition, $\bar{u}\neq 0$, $\overline {uv}=\bar{u}\bar{v}$ and $\mathrm{deg} (uv)=\mathrm{deg} (u)+\mathrm {deg} (v)$ for $u, v\neq 0$, as $\U$ is a domain.

The hypotheses of~\cite[Prop.\ 3.2]{LnLe05pr} can be slightly  weakened to yield, with essentially the same proof, the following  proposition.

\begin{proposition}\label{P:aut:ngr}
Let $A=\bigoplus_{i\in\N}A_{i}$ be an $\N$-graded $\K$-algebra with  $A_{0}=\K$ which is generated as an algebra by $A_{1}=\K x_{1}\oplus \cdots\oplus\K x_{n}$. Assume that for each $i\in\{ 1, \ldots , n \}$  there exist $0\neq a\in A$ and a scalar $q_{i, a}\neq 1$ such that
$x_{i}a=q_{i, a}ax_{i}$.
Then, given an algebra automorphism $\sigma$ of $A$ and a nonzero  homogeneous element $x$ of degree $d$, there exist $y_{d}\in A_{d} \setminus\{ 0 \}$ and $y_{>d}\in A_{\geq d+1}$ so that $\sigma (x)=y_ {d}+y_{>d}$.
\end{proposition}

The algebra $\U$, endowed with the grading just defined, satisfies  the conditions of the above proposition. Indeed, the quantum  Serre  relations involving $i$ and $i+1$ are equivalent to
\begin{align}
&e_{i}\left( e_{i}e_{i+1}-q^{-1}e_{i+1}e_{i} \right)=q\left( e_{i}e_{i +1}-q^{-1}e_{i+1}e_{i} \right)e_{i}\\
&e_{i+1}\left( e_{i}e_{i+1}-q^{-1}e_{i+1}e_{i} \right)=q^{-1}\left( e_ {i}e_{i+1}-q^{-1}e_{i+1}e_{i} \right)e_{i+1}.
\end{align}
Thus we have an analogue of~\cite[Cor.\ 3.3]{LnLe05pr}:

\begin{corollary}\label{C:aut:ngr}
Let $\sigma\in\aut$ and $x\in U_{d}\setminus\{ 0 \}$. Then $\sigma (x) =y_{d}+y_{>d}$, for some
$y_{d}\in U_{d}\setminus\{ 0 \}$ and $y_{>d}\in U_{\geq d+1}$.
\end{corollary}

\subsection{Invariance of the normal elements}\label{S:aut:inv}

\begin{proposition}\label{P:aut:inv:di}
Given $\sigma\in\aut$, there exist $\epsilon\in\{ 0, 1 \}$ and  nonzero scalars $\mu_{1}$ and $\mu_{3}$ such that $\eta^{\epsilon} \circ\sigma (\Delta_{i})=\mu_{i}\Delta_{i}$ for $i=1, 3$.
\end{proposition}
\begin{proof}
Since $\Delta_{1}$ is normal, so is $\sigma (\Delta_{1})$. By Lemma~\ref{L:int:nor} there exist $i\in\{ 1, 3\}$, $c\in\N$ and a central  element $z$ such that $\sigma (\Delta_{1})=\Delta_{i}^{c}z$.
Furthermore, $c\geq 1$ as $\Delta_{1}$ is not central. It follows  from Corollary~\ref{C:aut:ngr} that $c=1$, as $\mathrm{deg}(\Delta_ {j})=3$ for $j=1, 3$. Thus,
\begin{equation}\label{E:P:aut:inv:di}\sigma (\Delta_{1})=\Delta_{i}z. \end{equation}
If we repeat the argument above replacing $\Delta_1$ by $\Delta_i$ and
$\sigma$ by its inverse, apply $\sigma^{-1}$ to equation~(\ref {E:P:aut:inv:di}) and compute degrees, we find that $z$ is a  (nonzero) scalar. This same result can be reached by noticing that $ \Delta_{1}$ generates a (completely) prime ideal of $\U$, and so the  normal element $\sigma (\Delta_{1})$ must also generate such an  ideal. This, as well, implies that $z\in\K^{*}$. Similarly, $\sigma  (\Delta_{3})$ is a nonzero scalar multiple of $\Delta_{j}$ for some $j \in\{ 1, 3\}$ with $j\neq i$. If $i=1$ and $j=3$, we take $\epsilon=0 $; if $i=3$ and $j=1$, we take $\epsilon=1$. In either case, as $\eta$  interchanges $\Delta_{1}$ and $\Delta_{3}$, $\eta^{\epsilon}\circ \sigma$ fixes $\Delta_{1}$ and $\Delta_{3}$ up to scalars.
\end{proof}

We have as a corollary of Proposition~\ref{P:aut:inv:di} that any
algebra automorphism of $\U$ acts on the central element $z_{1}=
\Delta_{1}\Delta_{3}$ as multiplication by a scalar. Since the center
of $\U$ is $\K[z_{1}, z_{2}]$ with $z_{2}=\Delta_{2}$ and any $\sigma
\in\aut$ induces an automorphism of this polynomial algebra, it is
not hard to see that $\sigma (\Delta_{2})=\lambda \Delta_{2}+p
(z_{1})$ with $\lambda \in \K^*$ and $p(z_{1})$ a polynomial in $z_{1}$ with zero constant  term (by Corollary~\ref{C:aut:ngr}). Unfortunately, this is not quite  sufficient. In fact, if -- as we claim -- $\aut$ is the semi-direct  product of $\mathcal{H}$ and the order $2$ group generated by $\eta$,  it must be that $p(z_{1})=0$. Our next result, preceded by a  preparatory lemma, provides this step.

\begin{lemma}\label{L:aut:inv:ac}
For any $\sigma\in\aut$ there exist $\epsilon\in\{ 0, 1 \}$ and $\bar {\lambda}\in\mathcal{H}$ such that
\begin{equation}\left(\phi_{\bar{\lambda}}\circ\eta^{\epsilon}\circ \sigma-Id\right)\left( U_{1}\right)\subseteq U_{\geq 2}.\end{equation}
\end{lemma}
\begin{proof}
By Proposition~\ref{P:aut:inv:di}, $\eta^{\epsilon}\circ\sigma(\Delta_ {1})=t\Delta_{1}$, for some $\epsilon\in\{ 0, 1 \}$ and $t\in\K^{*}$.  Let $\psi=\eta^{\epsilon}\circ\sigma$. By Corollary~\ref{C:aut:ngr},  there exist $u_{1}\in U_{1}\setminus\{ 0\}$ and $u_{>1}\in U_{\geq 2} $ such that $\psi (e_{1})=u_{1}+u_{>1}$. If now we  apply $\psi$ to  the relation $e_{1}\Delta_{1}=q\Delta_{1}e_{1}$ and equate the  homogeneous terms of degree $4$, we obtain $u_{1}\Delta_{1}=q\Delta_{1} u_{1}$. As $u_{1}$ is a linear combination of $e_{1}$, $e_{2}$ and $e_ {3}$, Theorem~\ref{T:int:nor}(a) implies that $u_{1}=\lambda_{1}e_{1} $ for some $\lambda_{1}\in\K^{*}$. Analogously, $\psi (e_{i})=\lambda_ {i}e_{i}+w_{i}$ for $\lambda_{i}\in\K^{*}$ and $w_{i}\in U_{\geq 2}$,  $i=2, 3$. Let $\bar{\lambda}=(\lambda^{-1}_{1}, \lambda^{-1}_{2},  \lambda^{-1}_{3})$. Then $\left(\phi_{\bar{\lambda}}\circ\psi-Id \right)\left( U_{1}\right)\subseteq U_{\geq 2}$, since $\phi_{\bar {\lambda}}\left( U_{\
 !
 geq 2}\right)\subseteq U_{\geq 2}$.
\end{proof}

\begin{theorem}\label{T:aut:inv:d2}
Let $\sigma$ be an algebra automorphism of $\U$. Then there is a  nonzero scalar $\mu_{2}\in\K^{*}$ such that $\sigma(\Delta_{2})=\mu_ {2}\Delta_{2}$.
\end{theorem}
\begin{proof}
Since the statement of the theorem is valid for the automorphisms $ \eta$ and $\phi_{\bar{\lambda}}$, $\bar{\lambda}\in\mathcal{H}$, we  can assume by the previous lemma that $\left(\sigma-Id\right)\left( U_ {1}\right)\subseteq U_{\geq 2}$. Thus, by~\cite[Lem.\ 1.4.2]{AC92},  there exist $\mathrm{d}_{l}\in\mathrm{D}(\U)$, $l\geq 0$, such that
\begin{equation}\label{E:T:aut:inv:d2}\sigma (\Delta_{2})=\sum_{l\geq  0}\mathrm{d}_{l}(\Delta_{2}),\end{equation}
where $\mathrm{D}(\U)$ is the $\K$-subalgebra of $\mathrm{End}_{\K} \left(\U\right)$ generated by the $\K$-derivations of $\U$.  Furthermore, $\mathrm{d}_{0}(\Delta_{2})=\Delta_{2}$ and $\mathrm{d}_ {l}(\Delta_{2})$ is the homogeneous component of $\sigma (\Delta_{2}) $ of degree $l+4$, as $\Delta_{2}$ is homogeneous of degree $4$.

In Section~\ref{S:d} it will be shown (see Theorem~\ref
{T:der:der:der:der:der:main1}) that $\delta (\Delta_{2})$ is in the
ideal of $\U$ generated by $\Delta_{2}$, for any derivation $\delta$
of $\U$, and this will be done independently of Theorem~\ref
{T:aut:inv:d2}. Therefore, $\mathrm{d}(\Delta_{2})\in \left( \Delta_
{2}\right)$ for all $\mathrm{d}\in\mathrm{D}(\U)$ and thus $\sigma
(\Delta_{2})\in \left( \Delta_{2}\right)$, by~(\ref{E:T:aut:inv:d2}). 
This same reasoning applies to $\sigma^{-1}$, so that $\left(
\sigma (\Delta_{2}) \right)=  \left( \Delta_{2}\right)$. 
Since $\Delta_2$ is central, it is then obvious that there exists a
unit $\mu_2 \in \U$ such that $\sigma(\Delta_2)=\mu_2 \Delta_2$. 
However, the set of units of $\U$ is precisely $\K^*$, so that $\mu_2
\in \K^*$, as desired.
\end{proof}

\subsection{Determination of $\aut$}\label{S:aut:aut}

We are now ready to compute the group of algebra automorphisms of $\U$.

\begin{proposition}\label{P:aut:aut}
Let $\psi$ be an algebra automorphism of $\U$ with the property that $ \left(\psi-Id\right)\left( U_{1}\right)\subseteq U_{\geq 2}$. Then $ \psi$ is the identity morphism.
\end{proposition}
\begin{proof}
By the hypothesis on $\psi$, there exist $u_{i}\in U_{\geq (\mathrm {deg}(X_{i})+1)}$ such that
$$
\psi (X_{i})=X_{i}+u_{i}
$$
for all $1\leq i\leq 6$.
Also, by~Proposition~\ref{P:aut:inv:di} and Theorem~\ref {T:aut:inv:d2}, we know that $\psi (\Delta_{j})=\Delta_{j}$ for $j=1,  2, 3$. In particular, $u_{3}=0$ as $\Delta_{1}=X_{3}$.
Define, for $1\leq i\leq 6$, $d_{i}=\mathrm{deg}(\psi(X_{i}))$. It is
enough to prove that $d_{1}=d_{4}=d_{6}=1$ as $X_1=e_1$, $X_4=e_2$ and
$X_6=e_3$ generate $\U$ as an algebra. Let us assume, by way of  contradiction, that this is not the case. Thus $d_{1}+d_{4}+d_{6}>3$.

Notice that by Corollary~\ref{C:aut:ngr}, $d_{i}\geq\mathrm{deg}(X_ {i})$ for all $i$. Looking at the expression~(\ref{E:int:nor:d2}) of $ \Delta_{2}$ in the PBW basis and using the fact that $\psi$ fixes $ \Delta_{2}$, we can conclude that
\begin{equation}\label{E:aut:aut:d2f}d_{2}+d_{5}=d_{3}+d_{4}=3+d_{4}. \end{equation}
Also, since $X_{2}$ is a linear combination of $X_{1}X_{4}$ and $X_{4} X_{1}$, we have $2\leq d_{2}\leq d_{1}+d_{4}$ and similarly $2\leq d_ {5}\leq d_{4}+d_{6}$. Therefore,
\begin{align}\label{E:aut:aut:1}
d_{1}+d_{4}+d_{6}&\geq \mathrm{max}\{ d_{2}+d_{6}, d_{1}+d_{5} \}\q\q \mbox{and}\\\label{E:aut:aut:2}
d_{1}+d_{4}+d_{6}&>3=d_{3}.
\end{align}
Since $\psi$ fixes the degree $3$ element $\Delta_{3}$, the  inequality in~(\ref{E:aut:aut:1}) cannot be strict, by~(\ref {E:int:nor:d3}). Hence either $d_{1}+d_{4}+d_{6}=d_{2}+d_{6}$ or $d_ {1}+d_{4}+d_{6}=d_{1}+d_{5}$. These cases are symmetric and we can  assume without loss of generality that $d_{1}+d_{4}+d_{6}=d_{2}+d_{6}$.
Thus, using~(\ref{E:aut:aut:d2f}), $d_{1}+d_{4}=d_{2}=3+d_{4}-d_{5}$  and $d_{1}+d_{5}=3$. Since $d_{1}\geq 1$ and $d_{5}\geq 2$, it must  be $d_{1}=1$ and $d_{5}=2$. In other words, $u_{1}=0=u_{5}$ and $\psi $ fixes $X_{1}$ and $X_{5}$.

Now we apply $\psi$ to the defining equation~(\ref{E:int:nor:d2}) of $ \Delta_{2}$ to obtain
\begin{equation}\label{E:aut:aut:f1}u_{2}X_{5}=qX_{3}u_{4};\end {equation}
similarly, the relation $X_{5}X_{4}=q^{-1}X_{4}X_{5}$ yields
\begin{equation}\label{E:aut:aut:f2}X_{5}u_{4}=q^{-1}u_{4}X_{5}\end {equation}
after applying $\psi$; finally, $\psi$ applied to equation~(\ref {E:int:nor:d3}) gives
\begin{equation}\label{E:aut:aut:f3}\qhat \left( X_{1}X_{4}u_{6}+X_{1} u_{4}X_{6}+X_{1}u_{4}u_{6} \right)=
q\left( X_{2}u_{6}+u_{2}X_{6}+u_{2}u_{6} \right).\end{equation}
By~(\ref{E:aut:aut:f1}), $u_{2}=0\iff u_{4}=0$ and if this occurs then
$\qhat X_{1}X_{4}u_{6}=q X_{2}u_{6}$, on account of~(\ref {E:aut:aut:f3}). If $u_{6}\neq 0$ the latter implies
$\qhat X_{1}X_{4}=q X_{2}$, which is false as the $X_{i}$ form a PBW  basis. Thus $u_{6}=0$ and $d_{1}+d_{4}+d_{6}=3$, contradicting our  assumption. Hence $u_{4}, u_{2}\neq 0$. Likewise, if $u_{6}=0$ then~ (\ref{E:aut:aut:f3}) implies $\qhat X_{1}u_{4}=q u_{2}$ and then
by~(\ref{E:aut:aut:f1}) followed by~(\ref{E:aut:aut:f2}) we get $ \qhat X_{1}X_{5}u_{4}=qX_{3}u_{4}$, which is again a contradiction as  $u_{4}\neq 0$. Hence $d_{2}=\mathrm{deg}(u_{2})\geq 3$,
$d_{4}=\mathrm{deg}(u_{4})\geq 2$ and $d_{6}=\mathrm{deg}(u_{6})\geq 2$.

To obtain the final contradiction, we just have to look at the  degrees occurring
in~(\ref{E:aut:aut:f3}). Indeed, $\mathrm{deg}(X_{1}X_{4}u_{6})=2+d_ {6}<1+d_{4}+d_{6}=\mathrm{deg}(X_{1}u_{4}u_{6})$;  similarly, $\mathrm {deg}(X_{1}u_{4}X_{6})<\mathrm{deg}(X_{1}u_{4}u_{6})$, $\mathrm{deg} (X_{2}u_{6})<\mathrm{deg}(u_{2}u_{6})$ and $\mathrm{deg}(u_{2}X_{6})< \mathrm{deg}(u_{2}u_{6})$. Therefore we must have $\mathrm{deg}(X_{1} u_{4}u_{6})=\mathrm{deg}(u_{2}u_{6})$ and, using the notation  introduced in section~\ref{S:aut:ngr},
\begin{equation}\label{E:aut:aut:f4}\qhat\, X_{1}\bar{u}_{4}\bar{u}_ {6}=q\,\bar{u}_{2}\bar{u}_{6},\end{equation}so that $\qhat\, X_{1}\bar {u}_{4}=q\,\bar{u}_{2}$. Multiplying this equation on the right by $X_ {5}$, using relations $\bar{u}_{2}X_{5}=qX_{3}\bar{u}_{4}$ and $\bar {u}_{4}X_{5}=qX_{5}\bar{u}_{4}$, arising from~(\ref{E:aut:aut:f1})  and~(\ref{E:aut:aut:f2}), respectively, we obtain the equality $\qhat \, X_{1}X_{5}\bar{u}_{4}=qX_{3}\bar{u}_{4}$, which leads to the  contradiction
$\qhat\, X_{1}X_{5}=qX_{3}$. The contradiction was derived from the  assumption that $d_{1}+d_{4}+d_{6}>3$. Consequently $d_{1}=d_{4}=d_{6} =1$ and $\psi$ is the identity on $\U$.
\end{proof}

At last, we  prove our main result of this section, which gives a  positive
answer to the conjecture of Andruskiewitsch and Dumas~\cite{ADpr} for  $\U$.

\begin{theorem}\label{T:aut:aut:mt}
$\aut$ is isomorphic to the semi-direct product of the $3$-torus $ \mathcal{H}$ and the group of order $2$ generated by the diagram  automorphism $\eta$ of $\U$.
\end{theorem}
\begin{proof}
Let $\sigma\in\aut$. By Lemma~\ref{L:aut:inv:ac} and Proposition~\ref {P:aut:aut} there exist $\epsilon\in\{ 0, 1 \}$ and $\bar{\lambda}\in \mathcal{H}$ such that $\phi_{\bar{\lambda}}\circ\eta^{\epsilon}\circ \sigma$ is the identity on $\U$. Thus,
\begin{equation}\sigma=\eta^{\epsilon}\circ\phi_{\bar{\mu}},\end {equation}
where $\bar{\mu}=\bar{\lambda}^{-1}$. Furthermore,
the above expression is  easily seen to be unique, so the theorem  follows from~(\ref{E:aut:phieta}).
\end{proof}

\section{Derivations of $\U$}\label{S:d}

The aim of this section is to describe the Lie algebra of
$\K$-derivations of $\U$. In particular, we show that the Hochschild
cohomology group of degree 1 of $\U$ is a free module of rank 3 over
the center of $\U$. Our method consists of using previous results of
Osborn and Passman, \cite{OP95}, on the  the Hochschild
cohomology group of degree 1 of a quantum torus, and then to use the
theory of deleting derivations of Cauchon (see \cite{gC03}) in order
to transfer information on the derivations of a certain quantum torus
(in which $\U$ embeds) to  the derivations of $\U$ itself. This
method was first used in \cite{LL06} in order to describe the
derivations of the algebra of quantum matrices and of some related  algebras.

\subsection{The deleting derivations algorithm in $\U$}\label{S:der}

It follows from Section \ref{SS:int:rtp} that the theory of deleting
derivations (see \cite{gC03}) can be
applied to the iterated Ore extension
$R:=\U=\K[X_1]\dots[X_6;\tau_6,\delta_6]$. The corresponding
deleting derivations algorithm constructs, for each $r \in \{6,5,4,3,2 \}$, a family $(X_i^{(r)})_{i \in
\{1, \dots, 6 \} }$ of elements of $\fract(\U)$, defined as
follows (see \cite[Sec. 3.2]{gC03}): \\$ $
\begin{enumerate}
\item $X_1^{(6)}=X_1$, $X_2^{(6)}=X_2-q\hat{q}^{-1}X_3X_6^{-1}$, $X_3^ {(6)}=X_3$,
  $X_4^{(6)}=X_4-q\hat{q}^{-1}X_5X_6^{-1}$, $X_5^{(6)}=X_5$ and
  $X_6^{(6)}=X_6$.

In order to simplify the notations, we set $Y_i:=X_i^{(6)}$ for all $i
\in \{1,\dots,6\}$.

\item $X_1^{(5)}=Y_1-q\hat{q}^{-1}Y_3Y_5^{-1}$, $X_2^{(5)}=Y_2- qY_3Y_4Y_5^{-1}$, $X_3^{(5)}=Y_3$,
  $X_4^{(5)}=Y_4$, $X_5^{(5)}=Y_5$ and $X_6^{(5)}=Y_6$.

In order to simplify the notations, we set $Z_i:=X_i^{(5)}$ for all $i
\in \{1,\dots,6\}$.

\item $X_1^{(4)}=Z_1-q\hat{q}^{-1}Z_2Z_4^{-1}$, $X_2^{(4)}=Z_2$, $X_3^ {(4)}=Z_3$,
  $X_4^{(4)}=Z_4$, $X_5^{(4)}=Z_5$ and $X_6^{(4)}=Z_6$.

In order to simplify the notations, we set $T_i:=X_i^{(4)}$ for all $i
\in \{1,\dots,6\}$.

\item For all $r \in \{ 2,3 \}$ and $i \in \{1, \dots ,6 \}$, $X_i^ {(r)}=T_i$.
\end{enumerate}

As in \cite{gC03}, for all $r \in \{6,5,4,3,2\}$, we denote by $R^ {(r)}$ the subalgebra
of $\fract(R)$ generated by the elements $X_i^{(r)}$ for $i \in
\{ 1, \dots, 6 \}$. Also, we denote by $\overline{R}$ the subalgebra of
$\fract(R)$ generated by the indeterminates obtained at the end of this
algorithm, that is, $\overline{R}=R^{(2)}$ is the subalgebra of $ \fract(R)$
generated by the $T_i$, for each $i \in \{ 1, \dots, 6 \}$. Finally,  by convention, we set $R^{(7)}:=R$.\\

Recall from \cite[Th\'e. 3.2.1]{gC03} that, for all $r \in \{6,5,4,3,2 \}$, $R^{(r)}$ can be presented as an iterated Ore extension over
$\K$, with the generators $X_i^{(r)}$ adjoined in lexicographic order.
Thus the ring $R^{(r)}$ is a Noetherian domain.
Observe in particular that we have (with some abuse of notation):
\begin{equation}R^{(6)}=\K[Y_{1}][Y_{2};\tau_{2}][Y_{3};\tau_{3}][Y_ {4};\tau_{4},
    \delta_{4}][Y_{5};\tau_{5}, \delta_{5}][Y_{6};\tau_{6}],\end {equation}
\begin{equation}R^{(5)}=\K[Z_{1}][Z_{2};\tau_{2}][Z_{3};\tau_{3}][Z_ {4};\tau_{4},
    \delta_{4}][Z_{5};\tau_{5}][Z_{6};\tau_{6}],\end{equation}
\begin{equation}\overline{R}=R^{(4)}=R^{(3)}=R^{(2)}=\K[T_{1}][T_{2}; \tau_{2}][T_{3};\tau_{3}][T_{4};\tau_{4}][T_{5};\tau_{5}][T_{6};\tau_ {6}].
\end{equation}

Let $N \in \mathbb{N}^*$ and let $\Lambda=(\Lambda_{i,j})$
be a multiplicatively antisymmetric $N\times N$
matrix over $\K^*$;
that is, $\Lambda_{i,i}=1$ and $\Lambda_{j,i}=\Lambda_{i,j}^{-1}$
for all $i,j \in \{ 1,\dots ,N \}$.
We denote by $\K_{\Lambda}[T_1,\dots,T_N]$
the corresponding quantum affine space; that is, the $\K$-algebra
generated by the $N$ indeterminates $T_1,\dots,T_N$ subject to the
relations $T_i T_j =\Lambda_{i,j} T_j T_i $ for all $i,j \in  \{ 1, \dots ,N \}$. Next, we denote by  $\qtor$ the quantum
torus associated to the quantum affine space $\K_{\Lambda}[T_1, \dots,T_N]$, which is
the localisation of $\K_{\Lambda}[T_1,\dots,T_N]$ with respect to the  multiplicative system
generated by the $T_i$. For $\gamma=(\gamma_{1}, \ldots, \gamma_{N}) \in\Z^{N}$, set $T^{\gamma}:=T_1^{\gamma_1}  \dots
T_N^{\gamma_N}$. Note that the monomials $\left( T^{\gamma} \right)_ {\gamma\in\Z^{N}}$ form a PBW basis of $\qtor$.

It follows from \cite[Prop. 3.2.1]{gC03} that $\overline{R}$ is a  quantum affine space over $\K$ in the
indeterminates $T_1,\dots, T_6$. We denote by $\qtor$ the corresponding
quantum torus. In the present case, the matrix that defines the quantum
affine space $\overline{R}$ is the following:
$$\Lambda= \left(
\begin{array}{cccccc}
1      & q      & q      & q^{-1} & q^{-1} & 1 \\
q^{-1} & 1      & q      & q      & 1      & q^{-1} \\
q^{-1} & q^{-1} & 1      & 1      & q      & q \\
q      & q^{-1} & 1      & 1      & q      & q^{-1} \\
q      & 1      & q^{-1} & q^{-1} & 1      & q  \\
1      & q      & q^{-1} & q      & q^{-1} & 1 \\
\end{array}
\right) $$

For all $r \in \{6,5,4,3,2\}$, we denote by $S_r$ the
multiplicative system generated by the indeterminates $T_i$ with $i  \geq r$.
Since $T_i=X_i^{(r)}$ for all $i \geq r$, $S_r$ is a multiplicative  system of
regular elements of $R^{(r)}$. Moreover,  the $T_i$ with $i \geq r$  are normal in $R^{(r)}$. Hence
$S_r$ is an Ore set in $R^{(r)}$ and one can form
the localisation:
$$A_r:= R^{(r)}S_r^{-1}.$$
Clearly, the family $\left( (X_1^{(r)})^{\gamma_1}
(X_2^{(r)})^{\gamma_2} \dots(X_6^{(r)})^{\gamma_6} \right)$, with
$\gamma_i \in \mathbb{N} $ if $i < r$ and $\gamma_i  \in \mathbb{Z} $  otherwise, is a PBW basis of $A_{r}$.
Further, recall from \cite[Th\'e. 3.2.1]{gC03} that $\Sigma_r:=\{T_{r} ^k \mid k \in \mathbb{N} \}$ is
an Ore set in both $R^{(r)}$ and $R^{(r+1)}$, and that
$$
R^{(r)}\Sigma_r^{-1}=R^{(r+1)}\Sigma_r^{-1}.$$
Hence we get the following result.

\begin{lemma}\label{L:der:der:tower}
For all $r \in \{6,5,4,3,2\}$, we have $A_{r}=A_{r+1}\Sigma_r^{-1}$
with the convention that $A_7:=R=\U$.
\end{lemma}

Now, observe that $T_1$ is a normal element in $A_2$, so that one can
form the Ore localisation $A_1:=A_2 \Sigma_1^{-1}$, where $\Sigma_1$
is the multiplicative system generated by $T_1$. Naturally, $A_1$
is the quantum torus associated to the quantum affine space
$\overline{R}$. Hence we also denote $A_1$ by $\qtor$, and
we deduce from Lemma \ref{L:der:der:tower} the following tower of  algebras:
\begin{eqnarray}
\label{E:der:der:tower}
A_7=R & \subset & A_6=A_7 \Sigma_6^{-1} \subset  A_5= A_6 \Sigma_5^{-1}
\subset A_4= A_5 \Sigma_4^{-1}\\
& \subset & A_3= A_4 \Sigma_3^{-1} \subset
A_2= A_3 \Sigma_2^{-1} \subset A_1:= \qtor .
\end{eqnarray}

\subsection{Action of the deleting derivations algorithm on the normal
  elements} \label{Section32}

Observe that the formulas expressing the $Y_i$ in terms of the $X_i$
can be rewritten in order to express the $X_i$ in terms of the
$Y_i$. In particular, one can easily check that:
\\$X_1=Y_1$, $X_2=Y_2+q\hat{q}^{-1}Y_3Y_6^{-1}$, $X_3=Y_3$,
  $X_4=Y_4+q\hat{q}^{-1}Y_5Y_6^{-1}$, $X_5=Y_5$ and $X_6=Y_6$.

In a similar manner, one can express the $Y_i$ in terms of the $Z_i$,
and the $Z_i$ in terms of the $T_i$. More precisely, we have:
\\$Y_1=Z_1+q\hat{q}^{-1}Z_3Z_5^{-1}$, $Y_2=Z_2+qZ_3Z_4Z_5^{-1}$,
$Y_3=Z_3$, $Y_4=Z_4$, $Y_5=Z_5$ and $Y_6=Z_6$
\\and
\\$Z_1=T_1+q\hat{q}^{-1}T_2T_4^{-1}$, $Z_2=T_2$, $Z_3=T_3$,
  $Z_4=T_4$, $Z_5=T_5$ and $Z_6=T_6$.\\

Using these formulas, one can express
the three normal elements $\Delta_1$, $\Delta_2$ and $\Delta_3$
defined in Section \ref{SS:int:nor} in terms of the $Y_i$, or in  terms of the $Z_i$, or
in terms of the $T_i$. Indeed, straightforward computations lead to
the following results.

\begin{lemma}\label{L:der:der:der:delta}
\begin{enumerate}
\item $\Delta_1=X_3=Y_3=Z_3=T_3$.
\item $\Delta_2= X_2X_5-qX_3X_4=Y_2Y_5-qY_3Y_4=Z_2Z_5=T_2T_5$.
\item \begin{align*}
\Delta_{3}&=\qhat{\,}^{2}X_{1}X_{4}X_{6}-q\qhat X_{2}X_{6} -q\qhat
X_{1}X_{5}+q^{2}X_{3}\\
&=\qhat{\,}^{2}Y_{1}Y_{4}Y_{6}-q\qhat Y_{2}Y_{6}\\
&=\qhat{\,}^{2}Z_{1}Z_{4}Z_{6}-q\qhat Z_{2}Z_{6}\\
&=\qhat{\,}^{2}T_{1}T_{4}T_{6}
\end{align*}
\end{enumerate}
\end{lemma}

\subsection{Centers of the algebras $A_i$}

First, recall that the center of
$\U=A_7$ has been computed by  Alev and Dumas~\cite{AD94} and by
Caldero~\cite{pC94a, pC95}, who have shown that this is the  polynomial algebra
$\K[z_1,z_2]$, where $z_1=\Delta_1 \Delta_3$ and $z_2=\Delta_2$.

On the other hand, the center of the quantum torus $A_1=\qtor$ is
easy to compute. Indeed, it is well known (see for instance \cite{GL98})
that it is a Laurent polynomial ring over $\K$, and that it is
generated by the monomials $T_1^{\gamma_1} T_2^{\gamma_2} \dots T_6^ {\gamma_6}$, with
$\gamma_i \in \mathbb{Z} $, that are central. Easy computations show
that such a monomial is central if and only if
$\gamma_1=\gamma_4=\gamma_6=\gamma_3$ and $\gamma_2 =
\gamma_5$. Hence, we deduce from Lemma \ref{L:der:der:der:delta} that
the center of $\qtor$ is the Laurent polynomial ring over $\K$
generated by $z_1$ and $z_2$, that is:
\begin{equation*}
Z(\qtor) =Z(A_1)=\K[z_1^{\pm 1},z_2^{\pm 1}] .
\end{equation*}
It will be convenient  to denote by $\mathcal{F}$
the set of all $\gamma \in \mathbb{Z}^6$ such that $T^{\gamma}
\in Z(\qtor)$, that is:
\begin{equation}\label{E:cent:ai:def:f}
\mathcal{F}= \{\gamma \in \mathbb{Z}^6 \mid \gamma_1=\gamma_4= \gamma_6=\gamma_3
\mbox{ and }\gamma_2 =\gamma_5 \}.
\end{equation}

In the sequel we will also need  to know the center of $A_4$. Recall  that
$A_4$ is the localisation of the quantum affine space $R^{(4)}= \overline{R}$ at the multiplicative system generated by $T_4$, $T_5$  and $T_6$. In particular, the monomials $(T_1^{\gamma_1} T_2^ {\gamma_2} \dots T_6^{\gamma_6})$, with
$\gamma_i \in \mathbb{N} $ if $i \leq 3$ and $\gamma_i \in \mathbb{Z}$
otherwise, form a linear basis of $A_4$. The argument  used above to  compute the center of $\qtor$ also  works for $A_4$, with the  additional restrictions that $\gamma_i\geq 0$ for $i\leq 3$. So we  have the following result.

\begin{lemma}\label{L:der:der:der:der:center}
\begin{enumerate}
\item $Z(A_4)=Z(A_{7})= \K[z_1,z_2]$.
\item $Z(A_1)=\K[z_1^{\pm 1},z_2^{\pm 1}]$.
\end{enumerate}
\end{lemma}

\subsection{Derivations of $\U$}

Our aim in this section is to investigate the Lie algebra of $\K$-derivations of
$\U$, which we denote by $\der(\U)$.

Let $D$ be a derivation of $\U=A_7$. It follows from Lemma
\ref{L:der:der:tower} that $D$ extends (uniquely)
to a derivation of each of the algebras in the tower
$$A_7 \subseteq A_6 \subseteq \dots \subseteq A_2 \subseteq A_1= \qtor
.$$
In particular, $D$ extends to a derivation of the quantum torus
$\qtor$. So it follows from \cite[Cor. 2.3]{OP95} that $D$ can be
written as
$$D=\mathrm{ad}_x + \theta,$$
where $x \in\qtor$ and, in the
terminology of \cite{OP95}, $\theta$ is a central derivation of $\qtor $, that is,  $\theta(T_i)=\mu_i T_i$ with $\mu_i \in Z(\qtor)=\K[z_1^ {\pm 1},z_2^{\pm 1}] $.

Since the monomials $(T^{\gamma})_{\gamma \in \mathbb{Z}^6}$ form a PBW
basis of $\qtor$, one can write:
$$x= \sum_{\gamma \in \mathcal{E}} c_{\gamma} T^{\gamma},$$
where $\mathcal{E}$ is a finite subset of $\mathbb{Z}^6$ and
$c_{\gamma} \in \K$. Moreover, since
$\mathrm{ad}_x=\mathrm{ad}_{x+z}$ for all $z \in Z(\qtor)$,
it can be assumed that no monomial $T^{\gamma}$, with $\gamma\in \mathcal{E}$, belongs to $Z(\qtor)$, i.e., one can assume that $ \mathcal{E}\cap\mathcal{F}=\emptyset$.
Furthermore, by Lemmas~\ref{L:der:der:der:delta} and
\ref{L:der:der:der:der:center} we can write, for each $i\in\{ 1,  \ldots, 6 \}$, $\mu_i$ as follows:
$$\mu_i= \sum_{\gamma \in \mathcal{F}} \mu_{i,\gamma} T^{\gamma},$$
where $\mu_{i,\gamma}\in\K$.

\begin{lemma}\label{L:der:der:der:der:der:1}
For all $i \in \{1,2,3,4 \}$, we have $x \in A_i$.
\end{lemma}
\begin{proof}
We prove this lemma by induction on $i$. The case $i=1$ is trivial.  Hence we
assume that $x\in A_{i-1}$ for some $2\leq i \leq 4$.

It follows that
$$x= \sum_{\gamma \in \mathcal{E}} c_{\gamma} T^{\gamma},$$
where $\mathcal{E}$ is a finite subset of $\{ \gamma \in \mathbb{Z}^ {6}  \mid
\gamma_{1} \geq 0, \dots, \gamma_{i-2} \geq 0\}$ with  $\mathcal{E}  \cap \mathcal{F} =
\emptyset$.
We need to prove that $\gamma_{i-1} \geq 0$.

Let $j \in \{  1,\dots ,6 \}$ with $j \neq i-1$. As we have  previously observed, $D$ extends uniquely to a derivation of
$A_i$. Hence, since $T_{j} \in A_i$, we must have $D(T_{j}) \in A_i$,
that is:
\begin{eqnarray}
\label{v1beta}
xT_{j}-T_{j}x+\mu_{j}T_{j}\in A_i.
\end{eqnarray}
We set
$$x_+:= \sum_{\gamma \in \mathcal{E}, \gamma_{i-1} \geq 0} c_{\gamma}
T^{\gamma},$$
and
\begin{eqnarray}
\label{equa1}
x_-:= \sum_{\gamma \in \mathcal{E}, \gamma_{i-1}<0} c_{\gamma}
T^{\gamma}.
\end{eqnarray}
We shall prove that $x_-=0$.

First, we deduce from (\ref{v1beta}) that
$$u:=x_{-}T_{j}-T_{j}x_{-}+\mu_{j}T_{j}\in  A_i.$$
Next, using the commutation relations between the $T_k$,  we get
\begin{eqnarray}
\label{ubasis}
u = \sum_{\gamma \in \mathcal{E}, \gamma_{i-1}<0} c'_{j,\gamma}c_ {\gamma}
T^{\gamma+\varepsilon_{j}}+\sum_{\gamma \in  \mathcal{F}}
\mu'_{j,\gamma}T^{\gamma+\varepsilon_{j}}
\end{eqnarray}
where $\varepsilon_j$ denotes the $j$-th element of the canonical
basis of $\mathbb{Z}^6$,
$\mu'_{j,\gamma}=q^{\bullet}\mu_{j,\gamma}$ for some integer
$\bullet$, and  $c'_{j,\gamma} \in \K$ is defined by
$$x_-T_{j}-T_{j}x_-=\sum_{\gamma \in \mathcal{E}, \gamma_{i-1}<0} c'_ {j,\gamma}c_{\gamma}
T^{\gamma+\varepsilon_{j}}.$$
Observe that since we assume that $\mathcal{E} \cap \mathcal{F} =
\emptyset$, we have:
$$\mbox{for all }\gamma \in \mathcal{E} \mbox{\ and all }\gamma' \in  \mathcal{F}, \
\gamma+\varepsilon_{j} \neq \gamma' + \varepsilon_{j}.$$
Hence, (\ref{ubasis}) gives the expression of $u$ in the PBW basis of
$\qtor$.

On the other hand, since $u$ belongs to $A_i$, we get that:
$$u = \sum_{\gamma \in \mathcal{E}'} x_{\gamma}T^{\gamma},$$
where $\mathcal{E}'$ is a finite subset of $\{ \gamma \in \mathbb{Z}^ {6}  \mid
\gamma_{1} \geq 0, \dots, \gamma_{i-1} \geq 0\}$.
Comparing the two expressions of $u$ in the PBW basis of $\qtor$ leads
to  $c'_{j,\gamma}c_{\gamma}=0$ for all
  $\gamma \in \mathcal{E}$ such that $\gamma_{i-1} < 0$, as $j\neq  i-1$.
Hence, we have
$$x_-T_{j}-T_{j}x_-= \sum_{\gamma \in \mathcal{E}, \gamma_{i-1}<0}
c'_{j,\gamma}c_{\gamma} T^{\gamma+\varepsilon_{j}}=0,$$
for all $j \neq i-1$. In other words,
$x_-$ commutes with those $T_{j} $ such that $j \neq i-1$.

Now, recall from Lemma \ref{L:der:der:der:delta} that $z_1=\Delta_1
\Delta_3= \hat{q}^2 T_1T_4T_6T_3$ and $z_2=\Delta_{2}=T_2T_5$ are  central in
$\qtor$, so that $x_-$ commutes with those $T_{j} $ such that $j \neq  i-1$, and with $T_1T_4T_6T_3$ and $T_2T_5$. Naturally this implies
that $x_-$ also commutes with $T_{i-1}$, so that  $x_- \in Z(\qtor)$.  Thus one can write $x_-$ as follows:
\begin{eqnarray}\label{equa2}
x_-= \sum_{\gamma \in \mathcal{F}} d_{\gamma} T^{\gamma}.
\end{eqnarray}
As $\mathcal{E} \cap \mathcal{F}= \emptyset$, it follows from (\ref {equa1})
and (\ref{equa2}) that $x_-=0$, so that $x=x_+ \in  A_i$, as desired.
\end{proof}

In particular, it follows from Lemma \ref{L:der:der:der:der:der:1}
that $x \in A_4$. Since the derivation $D$  of $\U$ extends to a  derivation of
$A_4$, we must have $D(T_{i}) \in A_4$ for all $i \in \{1, \dots, 6 \} $. Hence
$$D(T_{i})=xT_{i}-T_{i}x+\mu_iT_{i} \in A_4.$$
Since  $x \in A_4$, this implies that
$\mu_{i}T_{i} \in A_4$ for all $i \in \{1, \dots, 6 \}$.
On the other hand, recall that $\mu_{i}$ is central in $\qtor$ and  can be written as:
$$\mu_{i}=\sum_{\gamma \in \mathcal{F}}\mu_{i,\gamma} T^{\gamma},$$
where $\mathcal{F}$ is given by~(\ref{E:cent:ai:def:f}).
Hence we get
\begin{equation*}
\begin{split}
\mu_{i}T_{i}&=\sum_{\gamma \in  \mathcal{F}}\mu'_{i,\gamma}
T^{\gamma+\varepsilon_{i}}\\ &=\sum_{\gamma =(\gamma_1,\gamma_2)\in
  \mathbb{Z}^2} \mu'_{i,\gamma}
T_1^{\gamma_1 +\delta_{1i}}T_2^{\gamma_2 +\delta_{2i}}T_3^{\gamma_1 + \delta_{3i}}T_4^{\gamma_1 +\delta_{4i}}T_5^{\gamma_2 +\delta_{5i}}T_6^ {\gamma_1 +\delta_{6i}} \in A_4,
\end{split}
\end{equation*}
where $\mu'_{i,\gamma}=q^{\bullet} \mu_{i,\gamma}$ for some integer $ \bullet$.

Assume now that $i \neq 2$. Then, since the monomials $T^{\gamma}$,  with $\gamma \in \mathbb{N}^3 \times
  \mathbb{Z}^3$, form a PBW basis of $A_4$ , we get that
  $\mu'_{i,\gamma}=0$ if either $\gamma_1 < 0$ or $\gamma_2 <0$.
Hence $\mu_i$ can be written as follows:
$$\mu_i = \sum_{\gamma=(\gamma_1,\gamma_2) \in \mathbb{N}^2}
c_{i,\gamma}T_1^{\gamma_1}T_2^{\gamma_2}T_3^{\gamma_1}T_4^{\gamma_1} T_5^{\gamma_2}T_6^{\gamma_1}.$$
In other words, $\mu_i \in \K[z_1,z_2] \subseteq \U$ since $z_1=\Delta_1
\Delta_3= \hat{q}^2 T_1T_4T_6T_3$ and $z_2=\Delta_{2}=T_2T_5$ by Lemma \ref{L:der:der:der:delta}.

Finally, assume that $i=2$. One cannot yet prove that $\mu_2 \in
\U=A_7$. However, one can prove the following weaker result: $\mu_2
z_2 \in  \K[z_1,z_2] \subseteq \U$. Indeed, we already know that
$\mu_2T_2 \in A_4$. Hence, it follows from Lemma \ref {L:der:der:der:delta} that
$\mu_2 z_2= \mu_2 T_2 T_5 \in A_4$. Further, $\mu_2 z_2$ is central in
$\qtor \supset A_4$, so that $\mu_2 z_2 \in Z(A_4)= \K[z_1,z_2]$, as  desired.

To sum up, we have just proved the following result.

\begin{corollary}\label{C:der:der:der:der:der:2}
\begin{enumerate}
\item $\mu_2 z_2 \in Z(A_4)=\K[z_1,z_2]
\subseteq \U$.
\item For all $i \neq 2$, $\mu_i \in \K[z_1,z_2]
\subseteq \U$.
\end{enumerate}
\end{corollary}

We now have to deal with localisation at elements which are not  normal. We do this in
three steps.

First, recall  from Lemma \ref{L:der:der:tower} that $A_4=A_5\Sigma_4^ {-1}$, where $\Sigma_4$
is the multiplicative system generated by $T_4=Z_4$.
Recall also that the monomials $Z_1^{\gamma_1} \dots Z_6^{\gamma_6}$,
with $\gamma =(\gamma_1,\dots,\gamma_6) \in \mathbb{N}^4 \times
\mathbb{Z}^2$, form a PBW basis of $A_5$. Of course, this implies that
the monomials $Z_1^{\gamma_1} \dots Z_6^{\gamma_6}$,
with $\gamma  \in \mathbb{N}^3 \times
\mathbb{Z}^3$, form a PBW basis of $A_4$.
In order to simplify the notation  we set, as usual,
$$Z^{\gamma}:=Z_{1}^{\gamma_{1}}Z_{2}^{\gamma_{2}} \dots
Z_{6}^{\gamma_{6}}$$
for all $\gamma \in \mathbb{N}^3 \times
\mathbb{Z}^3$.

\begin{corollary}\label{C:der:der:der:der:der:3}
$\mu_2 Z_2 \in A_5$.
\end{corollary}
\begin{proof} We know that $\mu_2 z_2 \in Z(A_4)=Z(A_5)$, so that
$\mu_2 z_2 \in A_5$. Now the result follows from the facts that
  $z_2=Z_2Z_5$ (Lemma \ref{L:der:der:der:delta}) and that $Z_5$ is  invertible in $A_5$.
\end{proof}

We are now able to prove that $x \in A_5$.

\begin{lemma}\label{L:der:der:der:der:der:4}
\begin{enumerate}
\item $x \in A_5$.
\item $\mu_2=\mu_1 +\mu_4 \in Z_{q}(\slf)$, where $Z_{q}(\slf)$ still denotes the center  of $\U$.
\item $D(Z_i)=\mathrm{ad}_x(Z_i) + \mu_i Z_i$ for all $i \in \{1, \dots ,6 \}$.
\end{enumerate}
\end{lemma}
\begin{proof} We proceed in three steps.\\

$\bullet$ {\it Step 1: We prove that $x \in A_5$.} \\

It follows from Lemma \ref{L:der:der:der:der:der:1} that $x $ belongs to
$A_4$, so that $x$ can be written as follows:
$$x= \sum_{\gamma \in \mathcal{E}} c_{\gamma} Z^{\gamma},$$
where $\mathcal{E}\subseteq\mathbb{N}^3 \times \mathbb{Z}^3$. \\ We set
$$x_+:= \sum_{\gamma \in \mathcal{E}, \gamma_{4} \geq 0} c_{\gamma}
Z^{\gamma},$$
and
$$x_-:= \sum_{\gamma \in \mathcal{E}, \gamma_{4} < 0} c_{\gamma}
Z^{\gamma}.$$
Assume that $x_- \neq 0$.

We denote by $B$ the subalgebra of
$A_4$ generated by the
$Z_{j}$ with $j \neq 4$, $Z_{5}^{-1}$ and $Z_{6}^{-1}$. Since $Z_4$ $q $-commutes with $Z_5$ and $Z_6$ in $A_4$, it is easy to
check that $A_4$ is a free left $B$-module with basis
$(Z_{4}^{a})_{a \in \mathbb{Z}}$, so that one can
write:
  $$x_-= \sum_{a =a_0}^{-1}  b_{a}
Z_{4}^{a} $$
with $a_0 <0$, $b_{a} \in B$ and $b_{a_0} \neq 0$. (Observe that this  makes sense since we are assuming that $x_- \neq 0$.)

As $D$  extends to a derivation of
$A_5$,  we have  $D(Z_{1}) \in A_5$. Recalling from Section
\ref{Section32} that
$Z_1=T_1+q\qhat^{-1} T_2T_4^{-1}$, this leads to:
$$x_-Z_{1}-Z_{1}x_-+\mu_{1}Z_{1} +q\hat{q}^{-1}(\mu_{2}-\mu_{1}-\mu_ {4})Z_{2}Z_{4}^{-1} \in A_5.$$
Since $\mu_1 \in \U \subset A_5$ by Corollary
\ref{C:der:der:der:der:der:2} and $Z_1 \in A_5$, we get
\begin{eqnarray}
\label{E:der:der:der:der:der:eqZ}
x_-Z_{1}-Z_{1}x_- +q\hat{q}^{-1}(\mu_{2}-\mu_{1}-\mu_{4})Z_{2}Z_{4}^ {-1} \in A_5.
\end{eqnarray}
Then, multiplying this expression by $Z_4$ (on the right) yields
$$(x_-Z_{1}-Z_{1}x_-)Z_{4} +q\hat{q}^{-1}(\mu_{2}-\mu_{1}-\mu_{4})Z_ {2} \in A_5.$$
Since $\mu_1$ and $\mu_4$ belong to $\U \subset A_5$ and $\mu_2Z_2 \in
A_5$ by Corollary \ref{C:der:der:der:der:der:3}, this leads to
$$u:=(x_-Z_{1}-Z_{1}x_-)Z_{4}  \in A_5,$$
that is:
$$u=   \sum_{a =a_0}^{-1}  b_{a}
Z_{4}^{a}Z_1 Z_4 - \sum_{a =a_0}^{-1} Z_{1} b_{a}
Z_{4}^{a+1} \in A_5. $$
Now, an easy induction shows that
$$Z_{4}^{-k}Z_{1}=q^{-k} Z_{1}Z_{4}^{-k}+q[k] Z_{2}Z_{4}^{-k-1}$$
for every positive integer $k$. Hence we have
$$u=  \sum_{a =a_0}^{-1} \left( q^{a}
b_{a} Z_1-Z_{1} b_{a} \right) Z_{4}^{a+1} +   \sum_{a
  =a_0}^{-1} q[-a]  b_{a}
Z_2Z_{4}^{a} \in A_5. $$

Since  $A_5$ is a free left $B$-module with basis
$(Z_{4}^{a})_{a \in \mathbb{N}}$ and $u \in
A_5$, one can write
$$u= \sum_{a =0}^{k}  u_{a} Z_{4}^{a} $$
with $k \in \mathbb{N}$ and $u_{a} \in B$.
Comparison of these two expressions of $u$ in the basis of $A_4$ (viewed
as a left $B$-module) shows that we must have $b_{a_0}=0$, a  contradiction.
Hence, $x_-=0$ and $x=x_+ \in A_5$, as desired.\\

$\bullet$ {\it Step 2: We prove that $\mu_{2}=\mu_{1}+\mu_{4}$.}\\

Since $x_-=0$, we deduce from (\ref{E:der:der:der:der:der:eqZ}) that
$$(\mu_{2}-\mu_{1}-\mu_{4})Z_{2}Z_{4}^{-1} \in A_5,$$
that is
$$(\mu_{2}-\mu_{1}-\mu_{4})Z_{2} \in A_5 Z_4.$$
Mutliplying this by $Z_5$ on the right leads to
$$(\mu_{2}-\mu_{1}-\mu_{4})z_2 \in A_5 Z_4,$$
since $z_2=Z_2Z_5$ by Lemma \ref{L:der:der:der:delta} and $Z_4 Z_5 =  q^{-1} Z_5 Z_4$.
We set $z:=(\mu_{2}-\mu_{1}-\mu_{4})z_2$ and  $J:=A_5Z_4$, so that $z
\in J$.

It follows from Corollary \ref{C:der:der:der:der:der:2} that $\mu_1,  \mu_4 \in \K[z_1,z_2]$ and
$\mu_2z_2 \in \K[z_1,z_2]$. Hence $z \in \K[z_1,z_2]$. We need to
prove that $z=0$. Let us write
$$z=\sum_{i,j \in \mathbb{N}}a_{i,j}z_1^i z_2^j,$$
with $a_{i,j} \in \K$ equal to zero except for a finite number of
them. Since $z_2 = Z_2 Z_5$ and $z_1=q^{-1}\qhat{\,}^{2}Z_3
Z_{1}Z_{6}Z_{4}-q\qhat Z_3Z_{2}Z_{6}$ (see Lemma \ref {L:der:der:der:delta}), we get that
$z_1 +q\qhat Z_{3}Z_{2}Z_{6} \in J$. Then, using the fact that $z_1$
and $z_2$ are central, we obtain that
$$z-\sum_{i,j \in \mathbb{N}}q^{\bullet}(-q\qhat)^i
a_{i,j}Z_2^{i+j} Z_3^i Z_5^j Z_6^i \in J,$$
where $\bullet$ denotes, as usual, an integer. Since we have already
proved that $z \in J$, this forces
\begin{equation}\label{sec:3:4:eq:1}
\sum_{i,j \in \mathbb{N}}q^{\bullet}(-q\qhat)^i
a_{i,j}Z_2^{i+j} Z_3^i Z_5^j Z_6^i \in J.
\end{equation}
However, since $Z_4$ $q$-commutes with $Z_5$ and $Z_6$, every element
of $J$ can be written as
\begin{equation}\label{sec:3:4:eq:2}
\sum_{\substack{\gamma \in \mathbb{N}^4 \times \mathbb{Z}^2\\
\gamma_4 >0}} c_{\gamma} Z_1^{\gamma_1} \dots Z_6^{\gamma_6}
\end{equation}
in the PBW basis of $A_5$.
Identifying the two expressions (\ref{sec:3:4:eq:1}) and (\ref{sec:3:4:eq:2}) leads to $a_{i,j}=0$ for all $i,j$, so that
$z=0$. Thus we have proved that
$(\mu_{2}-\mu_{1}-\mu_{4})z_{2}=0$. Since $z_2 \neq 0$, we get
$ \mu_{2}=\mu_{1}+\mu_{4}$, as desired. Observe that, since $\mu_1$
and $\mu_4$ belong to $Z_{q}(\slf)$ by Corollary \ref{C:der:der:der:der:der:2}, this implies that
$\mu_2$ also belongs to $Z_{q}(\slf)$.\\

$\bullet$ {\it Step 3: We prove that $D(Z_i)=\mathrm{ad}_x(Z_i) +\mu_i Z_i$ for
  all $i \in \{1, \dots, 6 \}$.}\\

If $i > 1$, this is trivial since $Z_i=T_i$ and we already know that
$D(T_i)=\mathrm{ad}_x(T_i) +\mu_i T_i$.

Next, recall that $Z_1=T_1+q\qhat^{-1}T_2T_4^{-1}$. Hence, we have
$$D(Z_1) = \mathrm{ad}_x(Z_1) + \mu_1T_1 + q\hat{q}^{-1} (\mu_2-\mu_4)T_2 T_4^ {-1}.$$
Since $\mu_2=\mu_1+\mu_4$, this implies that
$$D(Z_1) = \mathrm{ad}_x(Z_1) + \mu_1T_1 + q\hat{q}^{-1} \mu_1T_2 T_4^{-1}=\mathrm{ad}_x (Z_1) +
\mu_1Z_1,$$
as desired.

\end{proof}

We are now able to prove that $D(z_2)$ belongs to the ideal of $\U$
generated by $z_2=\Delta_2$. This result is crucial in order to compute
the automorphism group of $\U$ (see Theorem \ref{T:aut:inv:d2}).

\begin{theorem}\label{T:der:der:der:der:der:main1}
Let $D \in \der(\U)$. Then there exists $z \in Z_{q}(\slf)$ such that
$D(z_2)=zz_2$.
\end{theorem}
\begin{proof}
Let $D \in \der(\U)$. Since $z_2=\Delta_2=Z_2Z_5 \in A_5$ by Lemma
\ref{L:der:der:der:delta}, we deduce from Lemma \ref {L:der:der:der:der:der:4}  that
$D(z_2) = \mathrm{ad}_x(z_2) + (\mu_2+\mu_5)z_2$ with $x \in A_5$ and $\mu_2,
\mu_5 \in Z_{q}(\slf)$. Now the result easily follows  from the centrality of $z_2$
in $A_5$.
\end{proof}

Having completed the proof of Theorem \ref{T:aut:inv:d2} and thus  described
the automorphism group
of $\U$, we proceed to obtain a
complete description of $\der(\U)$.

Using arguments similar  to those in the
proof of Lemma \ref{L:der:der:der:der:der:4}, one can prove the
following two results.

\begin{lemma}\label{L:der:der:der:der:der:5}
\begin{enumerate}
\item $x \in A_6$.
\item $\mu_3 = \mu_1 +\mu_5$.
\item $\mu_2 + \mu_5=\mu_3 +\mu_4 $.
\item $D(Y_i)=\mathrm{ad}_x(Y_i) + \mu_i Y_i$ for all $i \in \{1, \dots ,6 \}$.
\end{enumerate}
\end{lemma}
And also:
\begin{lemma}\label{L:der:der:der:der:der:6}
\begin{enumerate}
\item $x \in A_7=\U$.
\item $\mu_3=\mu_2 +\mu_6$.
\item $\mu_5=\mu_4 +\mu_6$.
\item $D(X_i)=\mathrm{ad}_x(X_i) + \mu_i X_i$ for all $i \in \{1, \dots ,6 \}$.
\end{enumerate}
\end{lemma}

It is easy to check that we can define three derivations $D_1$, $D_4$  and
$D_6$ of $\U$ by setting:
$$\begin{array}{lllll}
D_1(X_1)= X_1 & D_1(X_2)= X_2 & D_1(X_3)= X_3 &
\\
D_4(X_2)= X_2 & D_4(X_3)= X_3 & D_4(X_4)= X_4 & D_4(X_5)=
X_5 \\
D_6(X_3)= X_3 &  D_6(X_5)=
X_5 & D_6(X_6)= X_6  &
\end{array}$$
and $D_{i}(X_{j})=0$ otherwise.

Then it follows from Lemmas \ref{L:der:der:der:der:der:4},
\ref{L:der:der:der:der:der:5} and \ref{L:der:der:der:der:der:6} that  any derivation $D $ of $\U$ can be
written as follows:
$$D=\mathrm{ad}_x + \mu_1 D_1+\mu_4 D_4 +\mu_6 D_6,$$
with $x \in \U$ and $\mu_1,\mu_4,\mu_6 \in Z_{q}(\slf)$.

Recall that the Hochschild cohomology group in degree 1 of $\U$,
denoted by $\mathrm{HH}^1(\U)$, is defined by:
$$\mathrm{HH}^1(\U):= \der(\U)/ \mathrm{InnDer}(\U),$$
where $ \mathrm{InnDer}(\U):=\{\mathrm{ad}_x \mid x \in \U\}$ is the Lie
algebra of inner derivations of $\U$. It is well known that
$\mathrm{HH}^1(\U)$ is a module over
$\mathrm{HH}^0(\U):=Z_{q}(\slf)$. Our final result makes  this latter  structure precise.

\begin{theorem}
\begin{enumerate}
\item Every derivation $D $ of $\U$ can be uniquely written as follows:
$$D=\mathrm{ad}_x + \mu_1 D_1+\mu_4 D_4 +\mu_6 D_6,$$
with $\mathrm{ad}_x \in \mathrm{InnDer}(\U)$ and $\mu_1,\mu_4,\mu_6 \in Z_{q} (\slf)$.
\item $\mathrm{HH}^1(\U)$ is a free $Z_{q}(\slf)$-module of rank 3 with
  basis $(\overline{D_1},\overline{D_4},\overline{D_6})$.
\end{enumerate}
\end{theorem}
\begin{proof}
It just remains to prove that, if $x \in \U$ and $\mu_1,\mu_4,\mu_6 \in
Z_{q}(\slf)$ with  $\mathrm{ad}_x + \mu_1 D_1+\mu_4 D_4 +\mu_6 D_6=0$, then
$\mu_1=\mu_4=\mu_6=0$ and $\mathrm{ad}_x=0$. Set $\theta:= \mu_1 D_1+\mu_4 D_4
+\mu_6 D_6 $, so that  $\mathrm{ad}_x + \theta =0 $. Since $\theta$ is a
derivation of $\U$, $\theta$ uniquely  extends to a derivation
$\tilde{\theta}$ of the quantum torus $\qtor$. Naturally, we still  have $\mathrm{ad}_x
+\tilde{\theta}=0$. Futher, straightforward
computations show that
$$\begin{array}{lll}
\tilde{\theta}(T_1)= \mu_1 T_1 & \tilde{\theta}(T_2)= (\mu_1+\mu_4)  T_2 &
\tilde{\theta}(T_3)=(\mu_1+\mu_4+\mu_6) T_3 \\ \tilde{\theta}(T_4)=  \mu_4 T_4 &
\tilde{\theta}(T_5)= (\mu_4+\mu_6)T_5 & \tilde{\theta}(T_6)= \mu_6 T_6
\end{array}$$
Hence $\tilde{\theta} $ is a central derivation of $\qtor$, in the  terminology
of \cite{OP95}. Thus we deduce from \cite[Cor. 2.3]{OP95} that
$\mathrm{ad}_x=0=\theta $. Evaluating $\theta$ on $X_1$,
$X_4$ and $X_6$ leads to $\mu_1=\mu_4=\mu_6=0$, as desired.
\end{proof}


\noindent
St\'ephane Launois\\ School of Mathematics, University of Edinburgh,\\
James Clerk Maxwell Building, King's Buildings, Mayfield Road,\\
Edinburgh EH9 3JZ, Scotland\\
E-mail: stephane.launois@ed.ac.uk \\[10pt]
Samuel A. Lopes\\
Centro de Matem\'atica da Universidade do Porto,\\
Rua do Campo Alegre 687, 4169-007 Porto, Portugal\\
E-mail: slopes@fc.up.pt

\end{document}